\input amstex
\documentstyle{amsppt}
\magnification=\magstephalf \hsize = 6.5 truein \vsize = 9 truein
\vskip 3.5 in

\rightheadtext\nofrills {A TIGHTER BOUND FOR THE NUMBER OF WORDS OF MINIMUM LENGTH}
\NoBlackBoxes
\TagsAsMath

\def\label#1{\par%
        \hangafter 1%
        \hangindent .75 in%
        \noindent%
        \hbox to .75 in{#1\hfill}%
        \ignorespaces%
        }

\newskip\sectionskipamount
\sectionskipamount = 24pt plus 8pt minus 8pt
\def\sectionskip{\vskip\sectionskipamount}
\define\sectionbreak{%
        \par  \ifdim\lastskip<\sectionskipamount
        \removelastskip  \penalty-2000  \sectionskip  \fi}
\define\section#1{%
        \sectionbreak   
        \subheading{#1}%
        \bigskip
        }

\redefine\qed{{\unskip\nobreak\hfil\penalty50\hskip2em\vadjust{}\nobreak\hfil
    $\square$\parfillskip=0pt\finalhyphendemerits=0\par}}

        
        \let    \< = \langle
        \let    \> = \rangle

\define\op#1{\operatorname{\fam=0\tenrm{#1}}} 
        \define         \a              {\alpha}
        \redefine       \b              {\beta}
        \redefine       \d              {\delta}
        \redefine       \D              {\Delta}
        \define         \e              {\varepsilon}
        \define         \E              {\op {E}}
        \define         \g              {\gamma}
        \define         \G              {\Gamma}
        \redefine       \l              {\lambda}
        \redefine       \L              {\Lambda}
        \define         \n              {\nabla}
        \redefine       \var            {\varphi}
        \define         \s              {\sigma}
        \redefine       \Sig            {\Sigma}
        \redefine       \t              {\tau}
        \define         \th             {\theta}
        \redefine       \O              {\Omega}
        \redefine       \o              {\omega}
        \define         \z              {\zeta}
        \define         \k              {\kappa}
        \redefine       \i              {\infty}
        \define         \p              {\partial}
        \define         \vsfg           {\midspace{0.1 truein}}

\topmatter
\title A tighter bound for the number of words of minimum length in an automorphic orbit
\endtitle
\author Donghi Lee
\endauthor

\address {Department of Mathematics, Pusan National University, Jangjeon-Dong, Geumjung-Gu, Pusan 609-735, Korea}
\endaddress

\email {donghi\@pusan.ac.kr}
\endemail


\abstract {Let $u$ be a cyclic word in a free group $F_n$ of
finite rank $n$ that has the minimum length over all cyclic words
in its automorphic orbit, and let $N(u)$ be the cardinality of the
set $\{v: |v|=|u|$ and $v= \phi(u)$ for some $\phi \in \text
{Aut}F_n\}$. In this paper, we prove that $N(u)$ is bounded by a
polynomial function of degree $2n-3$ in $|u|$ under the hypothesis
that if two letters $x, \, y$ with $x \neq y^{\pm 1}$ occur in
$u$, then the total number of $x^{\pm 1}$ occurring in $u$ is not
equal to the total number of $y^{\pm 1}$ occurring in $u$. We also
prove that $2n-3$ is the sharp bound for the degree of polynomials
bounding $N(u)$. As a special case, we deal with $N(u)$ in $F_2$
under the same hypothesis.}
\endabstract
\endtopmatter


\document
\baselineskip=24pt

\heading 1. Introduction
\endheading

Let $F_n$ be the free group of a finite rank $n$ on the set
$\{x_1, x_2, \dots, x_n\}$. We denote by $\Sig$ the set of {\it
letters} of $F_n$, that is, $\Sig = \{x_1, x_2, \dots, x_n\}^{\pm
1}$. As in [1, 6], we define a {\it cyclic word} to be a
cyclically ordered set of letters with no pair of inverses
adjacent. The {\it length} $|w|$ of a cyclic word $w$ is the
number of elements in the cyclically ordered set. For a cyclic
word $w$ in $F_n$, we denote the automorphic orbit $\{\psi(w):
\psi \in \text {Aut}F_n\}$ by $\text {Orb}_{\text {Aut}F_n}(w)$.

The purpose of this paper is to present a partial solution of the
following conjecture proposed by Myasnikov--Shpilrain [7]:

\proclaim {Conjecture} Let $u$ be a cyclic word in $F_n$ which has
the minimum length over all cyclic words in its automorphic orbit
$\text {Orb}_{\text {Aut}F_n}(u)$, and let $N(u)$ be the
cardinality of the set $\{v \in \text {Orb}_{\text {Aut}F_n}(u):
|v|=|u|\}$. Then $N(u)$ is bounded by a polynomial function of
degree $2n-3$ in $|u|$.
\endproclaim

This conjecture was motivated by the complexity of Whitehead's
algorithm which decides whether, for given two elements in $F_n$,
there is an automorphism of $F_n$ that takes one element to the
other. Indeed, proving that $N(u)$ is bounded by a polynomial
function in $|u|$ would yield that Whitehead's algorithm
terminates in polynomial time with respect to the maximum length
of the two words in question (see [7, Proposition 3.1]).

Proposing this conjecture, Myasnikov--Shpilrain [7] proved that
$N(u)$ is bounded by a polynomial in $|u|$ in $F_2$. Later, Khan
[3] improved their result by showing that $N(u)$ has the sharp
bound of $8|u|-40$ for $|u| \ge 9$ in $F_2$, by which the
conjecture was settled in the affirmative for $F_2$. For a free
group of bigger rank, Kapovich--Schupp--Shpilrain [2] showed that
$N(u)$ is bounded by a constant depending only on $n$ for $u$
contained in an exponentially generic subset of $F_n$, and the
author [4] recently proved that $N(u)$ is bounded by a polynomial
function of degree $n(5n-7)/2$ in $|u|$ under the following

\proclaim {Hypothesis 1.1} (i) A cyclic word $u$ has the minimum
length over all cyclic words in its automorphic orbit $\text {\rm
Orb}_{\text {\rm Aut}F_n}(u)$.

(ii) If two letters $x_i$ (or $x_i^{-1}$) and $x_j$ (or
$x_j^{-1}$) with $i < j$ occur in $u$, then the total number of
$x_i^{\pm 1}$ occurring in $u$ is strictly less than the total
number of $x_j^{\pm 1}$ occurring in $u$.
\endproclaim

In the present paper, we prove under the same hypothesis that
$N(u)$ is bounded by a polynomial function of degree $2n-3$ in
$|u|$, and that $2n-3$ is the sharp bound for the degree of
polynomials bounding $N(u)$:

\proclaim {Theorem 1.2} Let $u$ be a cyclic word in $F_n$ that
satisfies Hypothesis 1.1. Then $N(u)$ is bounded by a polynomial
function of degree $2n-3$ in $|u|$.
\endproclaim

\proclaim {Theorem 1.3} Let $n \ge 2$ be arbitrary. Then there
exist a polynomial $p_n(t)$ of degree exactly $2n-3$ in $t$ and a
sequence $(u_l)$ of cyclic words in $F_n$ satisfying Hypothesis
~1.1 such that $|u_l| \rightarrow \infty$ as $l \rightarrow
\infty$ and such that $N(u_l) \ge p_n(|u_l|)$. Thus $2n-3$ is a
sharp bound for the degree of a polynomial in $|u|$ bounding
$N(u)$ from above, provided $u$ is a cyclic word in $F_n$ that
satisfies Hypothesis 1.1.
\endproclaim

As a special case, we deal with $N(u)$ in $F_2$:

\proclaim {Theorem 1.4} Let $u$ be a cyclic word in $F_2$ that
satisfies Hypothesis 1.1. Then $N(u) \le 8|u|-40$. Moreover there
exists a sequence $(u_l)$ of cyclic words in $F_2$ satisfying
Hypothesis ~1.1 such that $|u_l| \ge 9$, $|u_l| \rightarrow
\infty$ as $l \rightarrow \infty$ and such that
$N(u_l)=8|u_l|-40$. Thus $N(u)$ has the sharp bound of $8|u|-40$
for $|u| \ge 9$.
\endproclaim

The same technique as used in ~[4] is applied to the proofs of
these theorems. The proofs will appear in Sections ~3--5. In
Section ~2, we will establish a couple of technical lemmas which
play an important role in the proof of Theorem ~1.2.

Now we would like to recall several definitions. As in [4], a {\it
Whitehead automorphism} $\s$ of $F_n$ is defined to be an
automorphism of one of the following two types (cf. [5, 8]):

\roster \item"(W1)" $\s$ permutes elements in $\Sig$. \item"(W2)"
$\s$ is defined by a set $A \subset \Sig$ and a letter $a \in
\Sig$ with both $a, a^{-1} \notin A$ in such a way that if $x \in
\Sig$ then (a) $\s(x)=xa$ provided $x \in A$ and $x^{-1} \notin
A$; (b) $\s(x)=a^{-1}xa$ provided both $x,\, x^{-1} \in A$; (c)
$\s(x)=x$ provided both $x,\, x^{-1} \notin A$.
\endroster
If $\s$ is of type (W2), we write $\s=(A, a)$. By $(\bar A,
a^{-1})$, we mean a Whitehead automorphism $(\Sig - A - a^{\pm 1},
a^{-1})$. It is then easy to see that $(A, a)(w)=(\bar A,
a^{-1})(w)$ for any cyclic word $w$ in $F_n$.

We also recall the definition of the degree of a Whitehead
automorphism of the second type (see ~[4]):

\proclaim {Definition 1.5} Let $\s=(A, a)$ be a Whitehead
automorphism of $F_n$ of the second type. Put $A'=\{i: \text
{either} \ x_i \in A \ \text {or} \ x_i^{-1} \in A, \ \text {but
not both} \}$. Then the {\it degree of $\s$} is defined to be
$\max A'$. If $A'=\emptyset$, then the {\it degree of $\s$} is
defined to be zero.
\endproclaim

Let $w$ be a fixed cyclic word in $F_n$ that satisfies Hypothesis
1.1 (i). For two letters $x, \, y \in \Sig$, we say that {\it $x$
depends on $y$ with respect to $w$} if, for every Whitehead
automorphism $(A, a)$ of $F_n$ such that
$$a \notin \{x^{\pm 1}, y^{\pm 1}\}, \ \{y^{\pm 1}\} \cap  A \neq \emptyset,
\ \text{and} \ \exists v \in \text {Orb}_{\text {Aut}F_n}(w): |(A,
a)(v)|=|v|=|w|,$$ we have $\{x^{\pm 1}\} \subseteq A$. Then, as
shown in ~[4], if $x$ depends on $y$ with respect to $w$, then $y$
depends on $x$ with respect to $w$.

We then construct the {\it dependence graph $\G_w$ of $w$} as
follows: Take the vertex set as $\Sig$, and connect two distinct
vertices $x,\, y \in \Sig$ by a non-oriented edge if either
$y=x^{-1}$ or $y$ depends on $x$ with respect to ~$w$. Let $C_i$
be the connected component of $\G_w$ containing $x_i$. Clearly
there exists a unique factorization
$$w=v_1v_2\cdots v_t \ \text {(without cancellation)},$$
where each $v_i$ is a non-empty (non-cyclic) word consisting of
letters in $C_{j_i}$ with $C_{j_i} \neq C_{j_{i+1}}$ ($i$ mod
$t$). The subword $v_i$ is called a {\it $C_{j_i}$-syllable} of
$w$. By the {\it $C_k$-syllable length of $w$} denoted by
$|w|_{C_k}$, we mean the total number of $C_k$-syllables of $w$.
We also define $|w|_s$ as $|w|_s=\sum_{k=1}^n |w|_{C_k}$.

\noindent {\bf Example 1.6.} Consider the cyclic word
$u=x_1^2x_2^3 x_3^4 x_4^5$ in $F_4$. Letting $v=(\{x_2^{\pm 1}\},
x_1)(u)=x_1x_2^3x_1x_3^4 x_4^5$, $v$ is an automorphic image of
$u$ with $|v|=|u|$ (hence $\G_u=\G_v$). This implies that both
$x_3^{\pm 1}$ and $x_4^{\pm 1}$ do not depend on $x_2^{\pm1}$.
Also putting $v'=(\{x_2^{\pm 1}\}, x_3^{-1})(u)$, we have
$|v'|=|u|$, so that $x_1^{\pm 1}$ does not depend on $x_2^{\pm
1}$. Hence the connected component $C_2$ of $\G_u$ containing
$x_2$ consists of only $x_2^{\pm 1}$. This way we can show that
the dependence graph $\G_u=\G_v$ has four distinct connected
components, each $C_i$ of which contains only $x_i^{\pm 1}$. Thus
$|u|_{C_i}=1$ for each $1 \le i \le 4$ and so $|u|_s=4$, whereas
$|v|_{C_1}=2$, $|v|_{C_j}=1$ for each $2 \le j \le 4$ and so
$|v|_s=5$.

\noindent {\bf Example 1.7.} Consider the cyclic word
$u=x_1^2x_2^3x_3^2x_4x_3^{-1}x_4x_3x_4^3$ in $F_4$, of which the
dependence graph $\G_u$ has three distinct connected components
$C_1$, $C_2$, $C_3=C_4$. Putting $v=(\{x_2^{\pm 1}\},
x_3^{-1})^2(u)=x_1^2x_3^2x_2^3x_4x_3^{-1}x_4x_3x_4^3$, $v$ is an
automorphic image of $u$ with $|v|=|u|$, so $\G_u=\G_v$. While
$|u|_{C_i}=1$ for each $1 \le i \le 4$ and so $|u|_s=4$,
$|v|_{C_1}=|v|_{C_2}=1$, $|v|_{C_3}=|v|_{C_4}=2$ and so $|v|_s=6$.

\heading 2. Preliminary Lemmas
\endheading

Throughout this section, when we say that $\s=(A, a)$ is a
Whitehead automorphism of $F_n$ of degree $i$, the following
restriction is additionally imposed:
$$\text{$a=x_j^{\pm 1}$ with $j
> i$.}$$ For two automorphisms $\phi$ and
$\psi$ of $F_n$, by writing $\phi \equiv \psi$ we mean the
equality of $\phi$ and $\psi$ over all cyclic words in $F_n$, that
is, $\phi(v)=\psi(v)$ for any cyclic word $v$ in $F_n$. For a
cyclic word $v$ in $F_n$, we define $M_k(v)$, for $k=0, 1, \dots,
n-1$, to be the cardinality of the set $\Omega_k(v)=\{\phi (v):
\phi$ can be represented as a composition $\phi=\a_t \cdots \a_1$
($t \in \Bbb N$) of Whitehead automorphisms $\a_i$ of $F_n$ of the
second type such that $k=\deg \a_t \ge \deg \a_{t-1} \ge \cdots
\ge \deg \a_1$ and $|\a_i \cdots \a_1(v)|=|v|$ for all $i=1,
\dots, t \}$.

\proclaim {Lemma 2.1} Under the foregoing notation, $M_1(v)$ is
bounded by a polynomial function of degree $n-1$ in $|v|$.
\endproclaim

\demo {Proof} Let $\ell_i$ be the number of occurrences of
$x_i^{\pm 1}$ in $v$ for $i=1, \dots, n$. Clearly $$M_1(v) \le
M_1(x_1^2x_2^{\ell_ 2} \cdots x_{n-1}^{\ell_{n-1}}
x_n^{\ell_n+\ell_1-2}).$$ So it is enough to prove that
$M_1(x_1^2x_2^{\ell_ 2} \cdots x_{n-1}^{\ell_{n-1}}
x_n^{\ell_n+\ell_1-2})$ is bounded by a polynomial function in
$|v|$ of degree $n-1$. Noting that $|x_1^2x_2^{\ell_ 2} \cdots
x_{n-1}^{\ell_{n-1}} x_n^{\ell_n+\ell_1-2}|_s=n$, put
$$\text{$\Lambda=\{v': |v'|_s=n$ and $v' \in \Omega _0(x_1^2x_2^{\ell_ 2}
\cdots x_{n-1}^{\ell_{n-1}} x_n^{\ell_n+\ell_1-2})\}$.}$$
Obviously the cardinality of the set $\Lambda$ is $(n-1)!$.

Let $w \in \Omega_1(x_1^2x_2^{\ell_ 2} \cdots x_{n-1}^{\ell_{n-1}}
x_n^{\ell_n+\ell_1-2})$. Then for an appropriate $v' \in \Lambda$,
there exist Whitehead automorphisms $\s_i$ of degree $0$ and
$\t_j$ of degree $1$ such that
$$w=\t_q \cdots \t_1 \s_p \cdots \s_1 (v'),
\tag 2.1
$$
where $|\s_i \cdots \s_1(v')|=|v'|$ and $|\s_i \cdots \s_1(v')|_s
\ge |\s_{i-1} \cdots \s_1(v')|_s$ for all $1 \le i \le p$, and
\linebreak $|\t_j \cdots \t_1 \s_p \cdots \s_1(v')|=|v'|$ for all
$1 \le j \le q$. Here, the same reasoning as in [4, Lemma 4.1]
shows that $\s_i \s_{i'} \equiv \s_{i'} \s_i$ for all $1 \le i, i'
\le p$. Furthermore, the chain $\t_q \cdots \t_1$ in (2.1) can be
chosen so that, for $\t_{ij}=(A_{ij}, a_{ij})$,
$$\t_q \cdots \t_1=(\t_{r q_r} \cdots \t_{r 1}) \cdots (\t_{2 q_2} \cdots \t_{2 1}) (\t_{1 q_1} \cdots \t_{1 1}),
\tag 2.2
$$
where $A_{ij}= A_{ij'}$ for all $1 \le j, j' \le q_i$, and $x_1
\in A_{i 1} \subsetneq A_{i+1 1}$.

We may assume without loss of generality that the index $r$ in
(2.2) is minimum over all chains satisfying (2.1) and (2.2).
Clearly in (2.1)--(2.2) the element $v'$ in $\Lambda$, the
Whitehead automorphisms $\s_1, \dots, \s_p$, and the index $r$ are
determined by $w$; so we put
$$v'_w=v', \quad \psi_w=\s_p \cdots \s_1, \quad \text {and} \quad r_w=r.$$
It is easy to see that $r_w$ is at most $n-1$.

For $s=1, \dots, n-1$, put
$$\text{$L_s=$ the cardinality of the set $\{\psi_w (v'_w): w \in \Omega _1(x_1^2x_2^{\ell_ 2} \cdots
x_{n-1}^{\ell_{n-1}} x_n^{\ell_n+\ell_1-2})$ with $r_w=s\}$.}$$
Then in view of (2.1)--(2.2), we have
$$M_1(x_1^2x_2^{\ell_ 2} \cdots x_{n-1}^{\ell_{n-1}} x_n^{\ell_n+\ell_1-2})
\le 2^{(n-1)}|v|L_1 + 2^{2(n-1)}|v|^2L _2+ \cdots
+2^{(n-1)^2}|v|^{n-1}L_{n-1},$$ since the number of possible $
A_{ij}$'s and the indices $q_i$'s in (2.2) are less than or equal
to $2^{n-1}$ and $|v|$, respectively. Hence it is enough to prove
that each $L_s$ is bounded by a polynomial function in $|v|$ of
degree $n-s-1$. Due to the result of ~[4, Lemma 4.1], there is
nothing to prove for $s=1$. So let $s \ge 2$ and put $E_i=A_{i 1}
- A_{i-1 1}$ for $i=2, \dots, s$. This can possibly happen only
when $\psi_w=\s_p \cdots \s_1$ in (2.1) can be re-arranged so
that, for $\s_j=(B_j, b_j)$,
$$\psi_w=(\s_{t_{s+1}} \cdots \s_{t_s+1}) \cdots (\s_{t_2}\cdots \s_2)\s_1,
\tag 2.3
$$
where $b_1 \in \{x_1^{\pm 1}\}$, $b_j^{\pm 1} \in E_i$ and either
$B_j \subseteq E_i$ or $B_j \cap E_i = \emptyset$ provided
$t_{i-1} < j \le t_i$ ($t_1=1$), and $b_j^{\pm 1} \notin (\bigcup
\limits_{i=2}^{s} E_i +x_1^{\pm 1})$ and either $B_j \subseteq
(\bigcup \limits_{i=2}^{s} E_i +x_1^{\pm 1})$ or $B_j \cap
(\bigcup \limits_{i=2}^{s} E_i +x_1^{\pm 1}) = \emptyset$ provided
$t_s < j \le t_{s+1}$.

Now, for $i=2, \dots, s$, let
$$\text{$h_i$ be the half of the
cardinality of the set $E_i$.}$$ Put $h=\sum \limits_{i=2}^s h_i$.
It then follows from the result of [4, ~Lemma 4.1] that the number
of cyclic words obtained by $\s_{t_{j+1}} \cdots \s_{t_j+1}$
applied to $(\s_{t_{j}} \cdots \s_{t_{j-1}+1}) \cdots
(\s_{t_2}\cdots \s_2)\s_1(v_w')$ is bounded by $|v|^{h_{j+1}-1}$
provided $j=1, \dots, s-1$ and by $|v|^{n-(h+1)-1}$ provided
$j=s$. Moreover the number of cyclic words derived from $\s_1$
applied to $v_w'$ is bounded by $n-2$. Therefore we have from
(2.3) that
$$L_s \le (n-1)! \, (n-2) |v|^{h_2-1} \cdots |v|^{h_s-1} |v|^{n-h-2}=(n-1)! \, (n-2) |v|^{n-s-1},$$
which is a polynomial function in $|v|$ of degree $n-s-1$, as
required. \qed
\enddemo

\proclaim {Remark} The proof of Lemma ~2.1 can be applied without
further change if we replace consideration of a single cyclic word
$v$, the length $|v|$ of $v$, and the total number of occurrences
of $x_j^{\pm 1}$ in $v$ by consideration of a finite sequence
$(v_1, \dots, v_m)$ of cyclic words, the sum $\sum\limits_{i=1}^m
|v_i|$ of the lengths of $v_1, \dots, v_m$, and the total number
of occurrences of $x_j^{\pm 1}$ in $(v_1, \dots, v_m)$,
respectively.
\endproclaim

\proclaim {Lemma 2.2} Under the foregoing notation, for each $k=2,
\dots, n-1$, $M_k(v)$ is bounded by a polynomial function of
degree $n+k-2$ in $|v|$.
\endproclaim

\demo {Proof} Let $\ell_i$ be the number of occurrences of
$x_i^{\pm 1}$ in $v$ for $i=1, \dots, n$. Since
$$M_k(v) \le M_k(x_1^2 \cdots x_k^2 x_{k+1}^{\ell_{k+1}} \cdots x_{n-1}^{\ell_{n-1}} x_n^{\ell_n+\ell_1+ \cdots +\ell_k-2k}),$$
it suffices to show that $M_k(x_1^2 \cdots x_k^2
x_{k+1}^{\ell_{k+1}} \cdots x_{n-1}^{\ell_{n-1}}
x_n^{\ell_n+\ell_1+ \cdots +\ell_k-2k})$ is bounded by a
polynomial function in $|v|$ of degree $n+k-2$. As in the proof of
Lemma ~2.1, put $\Lambda=\{v': |v'|_s=n$ and $v' \in
\Omega_0(x_1^2 \cdots x_k^2 x_{k+1}^{\ell_{k+1}} \cdots
x_{n-1}^{\ell_{n-1}} x_n^{\ell_n+\ell_1+ \cdots +\ell_k-2k})\}$.

Let $w \in \Omega_k(x_1^2 \cdots x_k^2 x_{k+1}^{\ell_{k+1}} \cdots
x_{n-1}^{\ell_{n-1}} x_n^{\ell_n+\ell_1+ \cdots +\ell_k-2k})$.
Then for an appropriate $v' \in \Lambda$, there exist Whitehead
automorphisms $\g_i$ of $F_n$ such that
$$w=\g_q \cdots \g_{p+1} \g_p \cdots \g_1 (v'),
\tag 2.4
$$
where the length of $v'$ is constant throughout the chain on the
right-hand side, $\deg \g_i=0$ provided $1 \le i \le p$, $\deg
\g_i> 0$ provided $p < i \le q$, and $|\g_j \cdots \g_1(v')|_s \ge
|\g_{j-1} \cdots \g_1(v')|_s$ for all $1 \le j \le p$. Here, since
$\g_i \g_{i'} \equiv \g_{i'} \g_i$ for all $1 \le i, i' \le p$ by
the same reasoning as in [4, Lemma 4.1], we may assume that either
none of $\g_i$ for $1 \le i \le p$ has multiplier $x_1$ or
$x_1^{-1}$ or only $\g_1$ has multiplier $x_1$ or $x_1^{-1}$. So
(2.4) can be re-written as
$$w=\g_q \cdots \g_{p+1} \g_p \cdots \g_1 \g_0(v'),
$$
where $\g_0$ is either the identity or a Whitehead automorphism of
$F_n$ of degree $0$ with multiplier $x_1$ or $x_1^{-1}$, and none
of $\g_j$ for $1 \le j \le q$ has multiplier $x_1$ or $x_1^{-1}$.

Write
$$\g_0(v')=x_1u_1x_1u_2 \hskip 0.2in \text {without cancellation}.
\tag 2.5
$$
(Note that $u_1$ and $u_2$ are non-cyclic subwords in $\{x_2,
\dots, x_n\}^{\pm 1}$.) Let $F_{n+1}$ be the free group on the set
$\{x_1, \dots, x_{n+1}\}$. From (2.5) we construct a pair $(v_1,
v_2)$ of cyclic words $v_1, v_2$ in $F_{n+1}$ with
$|v_1|+|v_2|=2|v|$ as follows:
$$v_1=x_1 u_1 x_{n+1} u_1^{-1} \quad \text {and} \quad v_2=x_1 u_2 x_{n+1} u_2^{-1}.$$

For each $\g_j=( D_j, d_j)$ for $1 \le j \le q$, define a
Whitehead automorphism $\e_j$ of $F_{n+1}$ as follows:
$$\align
&\text{if $x_1^{\pm 1} \in D_j$, then $\e_j=(D_j +x_{n+1}^{\pm 1}, d_j)$;} \\
&\text{if only $x_1 \in  D_j$, then $\e_j=(D_j +x_1^{-1}, d_j)$;} \\
&\text{if only $x_1^{-1} \in D_j$, then $\e_j=(D_j -x_1^{-1} +x_{n+1}^{\pm 1}, d_j)$;} \\
&\text{if $x_1^{\pm 1} \notin D_j$, then $\e_j=(D_j, d_j)$.}
\endalign
$$
Then arguing as in the proof of [4, Lemma ~4.2], we have $|\e_j
\cdots \e_1(v_1)|+|\e_j \cdots \e_1(v_2)|=2|v|$ for all $1 \le j
\le q$. Moreover, by the construction of $\e_j$, $\e_j$ is a
Whitehead automorphism of $F_{n+1}$ of degree at most $k$, and the
defining set of $\e_j$ contains either both of $x_1^{\pm 1}$ or
none of $x_1^{\pm 1}$. This yields the same situation as for a
chain of Whitehead automorphisms of $F_{n+1}$ of maximum degree
~$k-1$. Hence by the induction hypothesis together with the Remark
after Lemma ~2.1, $M_k(x_1^2 \cdots x_k^2 x_{k+1}^{\ell_{k+1}}
\cdots x_{n-1}^{\ell_{n-1}} x_n^{\ell_n+\ell_1+ \cdots
+\ell_k-2k})$ is bounded by $(n-2)$ times a polynomial function in
$2|v|$ of degree $(n+1)+(k-1)-2=n+k-2$, as required. \qed
\enddemo

\heading 3. Proof of Theorem 1.2
\endheading

Without loss of generality we may assume that $u$ satisfies
further

\roster \item "(i)" The $C_n$-syllable length $|u|_{C_n}$ of $u$
is minimum over all cyclic words in the set $\{v \in \text {\rm
Orb}_{\text {\rm Aut}F_n}(u): |v|=|u|\}$.

\item "(ii)" If the index $j$ ($1 \le j \le n-1$) is such that
$C_j \neq C_k$ for all $k>j$, then the $C_j$-syllable length
$|u|_{C_j}$ of $u$ is minimum over all cyclic words in the set
$\{v \in \text {Orb}_{\text {Aut}F_n}(u): |v|=|u| \ \text{and} \
|v|_{C_k}=|u|_{C_k} \ \text{for all} \ k>j \}$.
\endroster
(Namely, we may assume that $u$ satisfies further the conditions
in [4, Hypothesis 1.3].) Let $u' \in \text {Orb}_{\text
{Aut}F_n}(u)$ be such that $|u'|=|u|$. Due to the result of [4,
Theorem ~1.4], there exist Whitehead automorphisms $\pi$ of the
first type and $\t_1, \dots, \t_s$ of the second type such that
$$u'=\pi \t_s \cdots \t_1 (u),$$
where $n-1 \ge \deg \t_s \ge \deg \t_{s-1} \ge \cdots \ge \deg
\t_1$, and $|\t_i \cdots \t_1 (u)|=|u|$ for all $i=1, \dots, s$.
This implies that
$$N(u) \le C(M_0(u)+M_1(u)+\cdots+M_{n-1}(u)),
\tag 3.1
$$
where $C$ is the number of Whitehead automorphisms of $F_n$ of the
first type (which depends only on $n$), and $M_k(u)$ is as defined
in Section ~2. The result of [4, Lemma 4.1] shows that $M_0(u)$ is
bounded by a polynomial function in $|u|$ of degree $n-2$. Also by
Lemmas 2.1 and 2.2, $M_k(u)$ for each $k=1, \dots, n-1$ is bounded
by a polynomial function in $|u|$ of degree $n+k-2$. Then the
required result follows from (3.1). \qed

\heading 4. Proof of Theorem 1.3
\endheading
In [7], Myasnikov--Shpilrain pointed out that experimental data
provided by C. Sims show that the maximum value of $N(u)$ in $F_3$
is $48|u|^3-480|u|^2+1140|u|-672$ if $|u| \ge 11$ and this maximum
value is attained at $u=x_1^2x_2^2x_3x_2^{-1}x_3x_2x_3^{\ell}$
with $\ell \ge 3$. Inspired by this observation, we let
$$u=x_1^2 x_2 (x_2 x_nx_2^{-1}x_n) x_2 x_3 (x_3 x_nx_3^{-1}x_n)^2 x_3 \cdots x_{n-1}(x_{n-1} x_nx_{n-1}^{-1}x_n)^{n-2} x_{n-1} x_n^{\ell}$$
with $\ell \gg 1$ in $F_n$. Note that $u$ satisfies Hypothesis
1.1. For this $u$, we will prove that $N(u)$ cannot be bounded by
a polynomial function in $|u|$ of degree less than $2n-3$. For
each $i=2, \dots, n-1$ and $j=1, \dots, n-1$, let
$$\s_i=(\{x_i^{\pm 1}, \dots, x_n^{\pm 1}\}, x_n^{-1}) \quad \text{and} \quad \t_j=(\{x_j, x_{j+1}^{\pm 1}, \dots, x_{n-1}^{\pm 1}\}, x_n^{-1});$$
then $\s_i$ and $\t_j$ are Whitehead automorphisms of $F_n$ of
degree $0$ and degree $j$, respectively. Then the total number of
cyclic words derived from automorphisms of $F_n$ of the form
$\t_{n-1}^{m_{n-1}} \cdots \t_1^{m_1}\s_{n-1}^{k_{n-1}}\cdots
\s_2^{k_2}$, where $k_i, m_j \le \frac {\ell}{2n-3}$, applied to
$u$ is $(\frac {\ell}{2n-3})^{2n-3}$. Hence $N(u)$ is at least
$(\frac {\ell}{2n-3})^{2n-3}$, which completes the proof. \qed

\heading 5. Proof of Theorem 1.4
\endheading
Let us assume that $u$ satisfies further \roster \item "(i)" The
$C_2$-syllable length $|u|_{C_2}$ of $u$ is minimum over all
cyclic words in the set $\{v \in \text {\rm Orb}_{\text {\rm
Aut}F_n}(u): |v|=|u|\}$.

\item "(ii)" If $C_1 \neq C_2$, then the $C_1$-syllable length
$|u|_{C_1}$ of $u$ is minimum over all cyclic words in the set
$\{v \in \text {Orb}_{\text {Aut}F_n}(u): |v|=|u| \ \text{and} \
|v|_{C_2}=|u|_{C_2}\}$.
\endroster
(Namely, assume that $u$ satisfies further the conditions in [4,
Hypothesis 1.3].) Note that $M_0(u)=1$ in $F_2$, where $M_0(u)$ is
as defined in Section ~2. Also every Whitehead automorphism of
$F_2$ of degree $1$ is equal to either $(\{x_1\}, x_2)$ or
$(\{x_1\}, x_2^{-1})$ over all cyclic words in $F_2$. Hence, in
view of [4, ~Theorem ~1.4], $N(u)$ is the same as the cardinality
of the set $\{v: v=\pi \t^k(u) \, (k \ge 0)$, where $\pi$ is a
permutation on $\Sig$ and $\t$ is either $(\{x_1\}, x_2)$ or
$(\{x_1\}, x_2^{-1})$ such that $|\t^i (u)|=|u|$ for all $i=1,
\dots, k\}$. Let
$$\text{$\Lambda(u)=\{v: v=\t^k(u) \, (k \ge 0)$, where
$\t$ is as above$\}$.}$$

Let $m$ be the number of occurrences of $x_1^{\pm 1}$ in $u$.
First consider the maximum value $N(u)$ over all $u$ with $m=2$.
If $m=2$, then $u$ is of the form either
$x_1x_2^{\ell_1}x_1^{-1}x_2^{\ell_2}$ or $x_1^2x_2^{\ell}$. Then
the cardinality of $\Lambda(x_1x_2^{\ell_1}x_1^{-1}x_2^{\ell_2})$
equals $1$ and that of $\Lambda(x_1^2x_2^{\ell})$ equals $|u|-1$.
Hence $N(u)$ has the maximum value at $u=x_1^2x_2^{\ell}$. For
$u=x_1^2x_2^{\ell}$ with $\ell \ge 3$, $N(u)=4(|u|-1)$, since
there are $8$ permutations on $\Sig$ and
$\t^j(x_1^2x_2^{\ell})=\pi \t^{\ell-j}(x_1^2x_2^{\ell})$ for $j
\ge \ell/2$, where $\t=(\{x_1\}, x_2^{-1})$ and $\pi$ is the
permutation that fixes $x_1$ and maps $x_2$ to $x_2^{-1}$.

Next consider the maximum value of $N(u)$ over all $u$ with $m=4$.
(Here note that if $m$ is odd, then any Whitehead automorphism of
degree $1$ cannot be applied to $u$ without increasing $|u|$;
hence the cardinality of $\Lambda(u)$ equals $1$.) It is not hard
to see that $\Lambda(u)$ has the maximum cardinality $|u|-5$ at
$u=x_1^2x_2x_1^{-1}x_2x_1x_2^{\ell}$. For
$u=x_1^2x_2x_1^{-1}x_2x_1x_2^{\ell}$ with $\ell \ge 3$,
$N(u)=8(|u|-5)$, since $8$ permutations on $\Sig$ applied to the
elements of $\Lambda(x_1^2x_2x_1^{-1}x_2x_1x_2^{\ell})$ induce all
different cyclic words. Obviously this is the maximum value of
$N(u)$ over all $u$ with $m=4$.

Finally note that the cardinality of $\Lambda(u)$ cannot be
greater than nor equal to $|u|-5$ for any $u$ with $m >4$. This
means that $N(u) < 8(|u|-5)$ for every $u$ with $m >4$. Therefore,
the maximum value of $N(u)$ over all $u$ is $8(|u|-5)$, which is
attained at $u=x_1^2x_2x_1^{-1}x_2x_1x_2^{\ell}$ with $\ell \ge
3$. \qed

\heading Acknowledgements
\endheading
The author is grateful to the referee for many helpful comments
and suggestions. This work was supported by Pusan National
University Research Grant, 2004.

\heading References
\endheading

\roster

\item"1." P. J. Higgins and R. C. Lyndon, Equivalence of elements
under automorphisms of a free group, {\it J. London Math. Soc.}
{\bf 8} (1974), 254--258.

\item"2." I. Kapovich, P. E. Schupp and V. Shpilrain, Generic
properties of Whitehead's Algorithm and isomorphism rigidity of
random one-relator groups, {\it Pacific J. Math.} {\bf 223}
(2006), 113--140.

\item"3." B. Khan, The structure of automorphic conjugacy in the
free group of rank two, Computational and experimental group
theory, 115--196, {\it Contemp. Math.}, 349, Amer. Math. Soc.,
Providence, RI, 2004.

\item"4." D. Lee, Counting words of minimum length in an
automorphic orbit, preprint;
\newline http://www.arxiv.org/math.GR/0311410

\item"5." R. C. Lyndon and P. E. Schupp, ``Combinatorial Group Theory'', Springer-Verlag, New York/Berlin, 1977.

\item"6." J. McCool, A presentation for the automorphism group of
a free group of finite rank, {\it J. London Math. Soc.} {\bf 8}
(1974), 259--266.

\item"7." A. G. Myasnikov and V. Shpilrain, Automorphic orbits in free groups, {\it J. Algebra} {\bf 269} (2003), 18--27.

\item"8." J. H. C. Whitehead, Equivalent sets of elements in a
free group, {\it Ann. of Math.} {\bf 37} (1936), 782--800.

\endroster
\enddocument